\documentclass[11pt]{article}
\usepackage{amsmath,amssymb,amsfonts,amsthm,graphicx}
\usepackage{hyperref}
\newtheorem{theorem}{Theorem}[section]
\newtheorem{fact}[theorem]{Fact}
\newtheorem{example}[theorem]{Example}
\newtheorem{fig}[theorem]{Figure}
\newtheorem{proposition}[theorem]{Proposition}

\newtheorem{lemma}[theorem]{Lemma}

\newtheorem{remark}[theorem]{Remark}
\newtheorem{definition}[theorem]{Definition}
\numberwithin{equation}{section}
\title{Dependence with complete connections and the Gauss-Kuzmin theorem for $N$-continued fractions}
\author{
    Dan Lascu\footnote{e-mail: lascudan@gmail.com.}
    \nonumber \\
        }
\begin{document}
\maketitle
\thispagestyle{empty}
\begin{abstract}
We consider a family $\{T_N:N \geq 1 \}$ of interval maps as generalizations of the Gauss transformation. For the continued fraction expansion arising from $T_N$, we solve its Gauss-Kuzmin-type problem by applying the theory of random systems with complete connections by Iosifescu.
\end{abstract}

\section{Introduction} \label{section1}
The purpose of this paper is to prove a Gauss–Kuzmin type problem for $N$-continued fraction expansions introduced by Burger et al. \cite{Burger-2008}.
In order to solve the problem, we apply the theory
of random systems with complete connections by Iosifescu \cite{IG-2009}.

Fix an integer $N \geq 1$.
In this paper, we consider a generalization of the Gauss transformation, i.e.,
\begin{equation}\label{1.1}
T_{N}(x):=
\left\{
\begin{array}{ll}
{\displaystyle \frac{N}{x}- \left\lfloor\frac{N}{x}\right\rfloor},&
{ x \in I:=[0, 1], x \neq 0,}\\
0,& x=0
\end{array}
\right.
\end{equation}
where $\left\lfloor \cdot \right\rfloor$ denotes the floor (or entire) function.

The generalized Gauss measure
$\displaystyle G_N (A) = \frac{1}{\log \{(N+1)/{N}\}} \int_{A} \frac{dx}{x+N}$,
$A \in {\mathcal{B}}_I = \sigma$-algebra of all Borel subsets of $[0, 1]$ is $T_N$-invariant, i.e., $G_N\left(T_N^{-1}(A)\right)=G_N(A)$ for any $A \in {\mathcal{B}}_I$. Define $a_1(x)=\left\lfloor\frac{N}{x}\right\rfloor$,
$x \in (0, 1]$, $a_1(0)=\infty$, and $a_n(x)=a_1\left(T_N^{n-1}(x)\right)$, $x \in I$, $n \in \mathbb{N}_+:=\{1,2, \ldots\}$, with $T_N^{0}(x)=x$.
By the very definitions, Burger et al. proved in \cite{Burger-2008} that any irrational  $0<x<1$ can be written in the form
\begin{equation}
x = \displaystyle \frac{N}{a_1+\displaystyle \frac{N}{a_2+\displaystyle \frac{N}{a_3+ \ddots}}} :=[a_1, a_2, a_3, \ldots]_N \label{1.2}
\end{equation}
where $a_n$'s are non-negative integers. We will call (\ref{1.2}) the \emph{$N$-continued fraction expansion of $x$}.
In \cite{DKW-2013}, Dajani et al. proved that $(I,{\mathcal B}_{I}, G_N, T_N)$ is an ergodic dynamical system.

The Perron-Frobenius operator of $T_N$ under a probability measure $\mu$ on ${\mathcal B}_{I}$ such that $\mu(T_N^{-1}(A)) = 0$ whenever $\mu(A) = 0$ is defined as the bounded linear operator $U$ on the Banach space $L^1(I, \mu):=\{f: I \rightarrow \mathbb{C} : \int_{I} |f |d\mu < \infty \}$ such that the following holds:
\begin{equation}
\int_{A}Uf \,d\mu = \int_{T_{N}^{-1}(A)}f\, d\mu \quad
\mbox{ for all }
A \in {\mathcal{B}}_{I},\, f \in L^1(I,\mu). \label{1.3}
\end{equation}
In particular, the Perron-Frobenius operator of $T_N$ under $G_N$ is
\begin{equation}
\{Uf\}(x) = \sum_{i \geq N}V_{N,i}(x)\,f\left(\frac{N}{x+i}\right), \quad
f \in L^1(I,G_{N}) \label{1.4}
\end{equation}
where $V_{N,i}:= \displaystyle \frac{x+N}{(x+i)\,(x+i+1)}$ for any $i \geq N$ and $n \in \mathbb{N}_+$ \cite{L-2016-1}.
An important fact is that for any probability measure $\mu$ on ${\mathcal{B}}_{I}$ such that $\mu \ll \lambda$, where $\lambda$ is the Lebesgue measure on ${\mathcal{B}}_{I}$, we have
\begin{equation}
\mu \left(T_{N}^{-n}(A)\right)
= \int_{A} \{U^nf\}(x) d G_{N}(x) \label{1.5}
\end{equation}
where $f(x):= \left(\log\left(\frac{N+1}{N}\right)\right) (x+N) h(x)$ for
$x \in I$ \cite{L-2016-1}.

The problem of finding the asymptotic of $T_N^{-n}(A)$ as $n \rightarrow \infty$ represents the Gauss-Kuzmin-type problem for $N$-continued fraction expansions.

\begin{theorem} \label{G-K}
\rm{(A Gauss-Kuzmin theorem for $T_N$)}
\textit{Fix an integer $N \geq 1$ and let
$(I,{\mathcal B}_I,T_N)$ be as above.
\begin{enumerate}
\item[(i)]
For a probability measure $\mu$ on $(I,{\mathcal B}_I)$,
let the assumption (A) as follows:
\[(A) \quad\mbox{ $\mu$ is non-atomic and has a Riemann-integrable density.}\]
Then for any probability measure $\mu$ which satisfies (A),
the following holds:
\begin{equation} \label{1.6}
\lim_{n \rightarrow \infty}\mu (G_N^n < x)
= \frac{1}{\log \{(N+1)/N\}}\log \frac{x+N}{N}, \quad x \in I.
\end{equation}
\item[(ii)]
In addition to assumptions of $\mu$ in (i),
if the density of $I \ni x\mapsto \mu([0,x])$ is  Lipschitz continuous,
then there exist two positive constants $q < 1$ and $\ell$ such that
for any $x \in I$ and $n \geq 1$, there exists $\theta$ with $|\theta|\leq \ell$,
the following holds:
\begin{equation} \label{1.7}
\mu \left(G^n_N < x\right) =
\frac{1 + \theta q^n }{\log \{(N+1)/N\}}\log \frac{x+N}{N}
\end{equation}
where $\theta:=\theta(N,\mu,n,x)$.
As a consequence, the $n$-th error term $e_n(N,\mu;x)$
of the Gauss-Kuzmin problem is obtained as follows:
\begin{equation} \label{1.8}
e_{n}(N,\mu;x)=
\frac{\theta q^n }{\log \{(N+1)/N\}}\log \frac{x+N}{N}.
\end{equation}
\end{enumerate}}
\end{theorem}

The paper is organized as follows.
In Section \ref{section2}, we recall random system
with complete connections (=RSCC), and show examples and properties.
In Section \ref{section3}, we prove Theorem \ref{G-K}. By using the ergodic behavior of the RSCC associated with $N$-continued
fraction expansion, we determine the limit of the sequence
$(\,\mu (G_N^n < x)\,)_{n \geq 1}$ of distributions as $n \rightarrow \infty$.

\section{Random systems with complete connections} \label{section2}
In this section, we introduce random system with complete connections (=RSCC)
and show its properties.

The first explicit formal definition of the concept of dependence with
complete connections was given by Onicescu and Mihoc \cite{OM-1935}.
It is a non-trivial extension of Markovian dependence theory,
and it was also investigated by Doeblin and Fortet \cite{DF-1937}
and by Harris \cite{Harris-1955}.
The concept of random system with complete connections
was defined by Iosifescu \cite{MI-1963}.

Examples of RSCC are urns models \cite{OM-1935,IT-1969},
stochastic learning processes \cite{Norman-1972, IT-1969, Karlin-1953},
partially observed random chains \cite{IT-1969}, image coding \cite{BarnsleyElton-1988},
continued fraction expansions \cite{IG-2009,L-2013,Sebe-2001,Sebe-2002,SebeLascu-2014} and others.

An RSCC is often called an
\textit{iterated function system with place-dependent probabilities}
or simply an \textit{iterated function system (=IFS)}
\cite{Barnsley-1988}.
Applications of IFS to continued fractions can be found
in \cite{KHL-2009,MU-2003}.
For more detail, see \cite{IG-2009, IT-1969, Kalpazidou-1986}.
\subsection{Definitions and their extensions} \label{subsection2.1}
In this subsection, we introduce the definition of random system with complete connections which is restricted to a smaller class than the original.
Next we extend domains of functions in such a system.
\begin{definition} \label{def2.1} \cite{MI-1963,IG-2009,Norman-1972}
A random system with complete connections (=RSCC) is a quadruple

\begin{equation}
\{(W, {\mathcal W}), X, u, P\} \label{2.1}
\end{equation}
where
\begin{enumerate}
\item[(i)] $(W, {\mathcal W})$ is a measurable space and $X$ is a non-empty set;
\item[(ii)] $u: W \times X \rightarrow W$ is a measurable function
with respect to $W$ for each element in $X$;
\item[(iii)] $P:W \times X \rightarrow [0, 1]$ satisfies
that $\sum_{x \in X} P(w,x)=1$ for each $w \in W$, and
$P(\cdot,x)$ is a measurable function on $(W,{\mathcal W})$ for each $x \in X$.
\end{enumerate}
\end{definition}

\noindent
For an RSCC in Definition \ref{def2.1}, we call $W$, $X$, $u$ and $P$,
the \textit{state space}, the \textit{event space}, the \textit{transition function}
and the \textit{($(W,X)$-) transition probability function}, respectively
(\cite{Norman-1972}, Section 1.2).
The role of the function $u$ is to change a state $w \in W$
into the new state $w^{'}=u(w,x)\in W$ by an event $x \in X$:
\begin{equation} \label{2.2}
W \ni w \stackrel{x}{\mapsto} w^{'}=u(w,x) \in W.
\end{equation}
In this case, $P(w,x)$ is regarded as the probability of the transition $w \mapsto w^{'}$
which depends on the information of both $w$ and $x$.
\begin{remark} \label{rem2.2}
{\rm
\begin{enumerate}
\item[(i)]
In the original definition of RSCC in \cite{IG-2009},
$X$ is assumed as a measurable space $(X,{\mathcal X})$.
In Definition \ref{def2.1}, we treat only the case ${\mathcal X}={\mathcal P}(X)$ (= the power set of $X$).
A lot of examples in \cite{IG-2009} satisfy this condition.
\item[(ii)]
An RSCC can be regarded as an automaton with output \cite{Mealy,Kawamura}.
A {\it Mealy machine} $(W,X,\Delta,u,\lambda)$ consists of
three sets $W,X,\Delta$ and two maps $u,\lambda$ where
$W$ is the \textit{set of states},
$X$ is the \textit{input alphabet},
$\Delta$ is the \textit{output alphabet}, and
the \textit{transition function} $u:W \times X^{*}\rightarrow W$ and
the \textit{map of outputs} $\lambda : W \times X^{*} \rightarrow \Delta^{*}$ which satisfy
\begin{equation} \label{2.3}
\left\{
\begin{array}{ll}
u(w,xa)=&u(u(w,x),a),\\
\\
\lambda(w,xa)=&\lambda(w,x)\lambda(u(w,x),a),
\end{array}
\right.
\quad x \in X^*,\,a \in X,w \in W
\end{equation}
where $X^{*}$ and $\Delta^{*}$ denote
free semigroups generated by $X$ and $\Delta$, respectively.
When it is an RSCC,
$\Delta=[0,1]$ and
$\Delta^{*}$ is reduced to the $[0,1]$ with respect to
the multiplication in $(\mathbb{R}, \times)$, $\pi:\Delta^{*}\rightarrow [0,1]$.
The map $\pi\circ \lambda$ is a transition probability function $P$.
An example of Mealy machine as an RSCC
will be shown in Example \ref{ex2.5}.
The similarity can be explained as the reason
that initial applications of RSCC were related to learning processes
which may be understood as an algorithmic study of such systems.
\end{enumerate}
}
\end{remark}

For a given RSCC $\{(W, {\mathcal W}), X, u, P\}$,
we extend domains of both $u$ and $P$ as follows.
\begin{enumerate}
\item[(i)]
In order to extend the domain of $u$, we prepare the free semigroup $X^{*}$
generated by $X$ as follows:
In Definition \ref{def2.1}, we regard $X$ as a set of letters,
and treat $X^{n}$ as the set of all strings of length $n \geq 1$.
We write an element $(x_1, \ldots, x_n) \in X^n$ as $x_1 \cdots x_n$.
Then the set $X^{*}$ of all strings with the empty letter $\varepsilon$
is a semigroup with the concatenation as its product,
and the unit of $X^{*}$ is $\varepsilon$.
The semigroup $X^{*}$ is called the \textit{free semigroup generated by }
$X$ \cite{Lothaire}.

\quad With respect to the semigroup $X^{*}$,
the function $u$ induces a unique right action of $X^{*}$ on $W$
which is denoted by the same symbol $u$ here:
\begin{equation} \label{2.4}
u:W \times X^{*}\to W; \quad (w,x) \mapsto wx:=u(w,x).
\end{equation}
For example,
$w(xx^{'})$ is defined as $(wx)x^{'}
=u(u(w,x),x^{'})$ for $w \in W$, $x \in X^*$ and $x^{'}\in X$.
For $x=x_{1}\cdots x_n \in X^n$,
we can write
$wx:=w(x_{1}\cdots x_n)$ for $w \in W$.
\item[(ii)]
Let ${\mathcal P}(X)$ denote the power set of $X$.
The domain $W \times X$ of $P$ is extended to $W \times {\mathcal P}(X)$
as follows:
\begin{equation} \label{2.5}
P(w,A):=
\left\{
\begin{array}{ll}
\displaystyle\sum_{x \in A}P(w,x), \quad &A \ne \emptyset,\\
\\
0, \quad &A=\emptyset,
\end{array}
\right.
\quad (w,A)\in W \times {\mathcal P}(X).
\end{equation}
By Definition \ref{def2.1}(iii),
we see that $(X,{\mathcal P}(X),P(w,\cdot))$ is a probability space for each $w \in W$.
\end{enumerate}

As a generalization of $P$,
we define $P_r:W \times X^{r} \to [0,1]$ for each $r \geq 1$ by
\begin{equation} \label{2.6}
P_1:=P,\quad
P_r(w,xx'):=P_{r-1}(w,x)P(wx,x'),
\end{equation}
with $w \in W,\,x\in X^{r-1},\,x \in X,\,r \geq 2$,
where the notation in (\ref{2.4}) is used.
Then we see that $P_r(w,x)$ means the transition probability
from $w$ to $w \cdot x$ along the path $w \rightarrow wx_1 \rightarrow wx_1x_2 \rightarrow \cdots \rightarrow w x_1 \cdots x_r$ with length $r$.
For $A \subset X^r$, we also define
\begin{equation} \label{2.7}
P_r(w, A) :=
\left\{
\begin{array}{ll}
\displaystyle \sum_{x \in A} P_r(w,x), \quad &A \ne\emptyset,\\
\\
0,\quad &A=\emptyset.
\end{array}
\right.
\end{equation}
By definition,
$P_r(w, A)$ is the summation of all paths from $w$ to $w \cdot x$
for $x \in A$.
For example, we can verify that
$P_r(w,A \times X)= P_{r-1}(w,A)$ for any
$(w,A)\in W \times {\mathcal P}(X^{r-1})$ by Definition \ref{def2.1}(iii) and (\ref{2.6}).
Especially,
$P_r(w,X^r)= 1$.
Hence
$(X^r,{\mathcal P}(X^r),P_r(w,\cdot))$ is also a probability space
for each $w \in W$ and $r \geq 1$.

\subsection{Transition probability functions defined on the square of the state space} \label{subsec2.2}
Let $\{(W,{\mathcal W}),X,u,P\}$ be as in
Definition \ref{def2.1}.
Recall that $P$ is the transition probability function
with the domain $W \times {\mathcal P}(X)$.
We define new transition probability functions
with the domain $W \times {\mathcal W}$ in this subsection.
\subsubsection{Space of paths} \label{subsubsec2.2.1}
For $w,w^{'} \in W$, define the subset $X(w,w^{'})$ of $X$ by
\begin{equation} \label{2.8}
X(w,w^{'}):=\{x \in X: u(w,x)=w^{'}\}.
\end{equation}
Then $X(w,w^{'})$ can be identified with the set of
all paths from $w$ to $w^{'}$ with length $1$.
Remark that $X(w,w^{'})$ may be the empty set.
By definition,
$X(w,w^{'})\cap X(w,w^{''})=\emptyset$
when $w^{'}\ne w^{''}$ and
$\coprod_{w^{'} \in W} X(w,w^{'})=X$.
From these properties, the following holds.
\begin{fact} \label{fact2.3}
For $(w,B) \in W \times {\mathcal W}$, let $X(w,B):=\{x \in X: u(w,x) \in B\}$.
Then the following holds for each $w \in W$:
\begin{enumerate}
\item[(i)]
$X(w,B)\cap X(w,B^{'})=\emptyset$ when $B \cap B^{'}=\emptyset$.
\item[(ii)]
$X(w,B)=\bigcup_{w^{'} \in B}X(w,w^{'})$.
\item[(iii)]
$X(w,B) \subset X(w,B^{'})$ when $B \subset B^{'}$.
\item[(iv)]
$X(w,B) \cup X(w,B^{'})=X(w,B \cup B^{'})$.
\end{enumerate}
\end{fact}

\noindent
In \cite{IG-2009}, the symbol $X(w,B)$ is written as $B_w$.

\subsubsection{Transition probability functions defined on the square of the state space} \label{subsubsec2.2.2}
By using the extension of $P$ in (\ref{2.5}) and
$X(w,w^{'})$ in (\ref{2.8}), define the new function $Q:W \times W \rightarrow [0,1]$ by

\begin{equation} \label{2.9}
Q(w, w^{'}): = P(w, X(w,w^{'})), \quad (w,w^{'}) \in W \times W.
\end{equation}
This means the total sum of transition probabilities
from $w$ to $w^{'}$ by all possible event $x \in X$.
We extend the domain of $Q$ to $W \times {\mathcal W}$
as follows (\cite{IG-2009}, (1.1.11)):
\begin{equation}\label{2.10}
Q(w, B): =P(w,X(w,B))= \sum_{x \in X(w,B)} P(w, x),
\quad (w,B) \in W\times {\mathcal W}.
\end{equation}
We call $Q$ the \textit{($(W,W)$-) transition probability function} associated with
$\{(W, {\mathcal W}), X, u, P\}$.
The number $Q(w, B)$ means the
probability such that a state $w$ changes to  a state belonging to $B$
by measuring all possible (one-step) events.
For $w \in W$, define the function $Q_w$ as
\begin{equation} \label{2.11}
Q_w:{\mathcal W} \to [0,1];\quad
Q_w(B):=Q(w,B),\quad B \in {\mathcal W}.
\end{equation}
By Fact \ref{fact2.3},
\begin{equation} \label{2.12}
0 \leq Q_w(B) \leq Q_w(W)=1,\quad
Q_w(B)+Q_w(B^{'})=Q_w(B \cup B^{'})\mbox{ when }B \cap B^{'}=\emptyset.
\end{equation}
Hence $(W,{\mathcal W},Q_w)$ is a probability space for each $w \in W$.
Remark that $Q_w$ can be defined on the whole of ${\mathcal P}(W)$,
but we restrict $Q_w$ on the given $\sigma$-algebra ${\mathcal W} \subset {\mathcal P}(W)$ here.
Thanks to $Q_w$, we can define the integration $\int_{B} f(w^{'})\, dQ_w(w^{'})$
for $B \in {\mathcal W}$ and a measurable function $f$ on $(W,{\mathcal W})$.

For $w,w^{'} \in W$, define the family $\{Q^{(k)}(w,w^{'}): k \geq 1\}$ recursively as
\begin{equation} \label{2.13}
Q^{(1)} (w, w^{'}):= Q(w, w^{'}),\quad
Q^{(k)} (w, w^{'}):=
\int_{W} dQ_w(w^{''})\,Q^{(k-1)}(w^{''}, w^{'}), \quad k \geq 2.
\end{equation}
By definition,
$Q^{(k)} (w, w^{'})$ is the expectation value of $Q^{(k-1)}(\cdot, w^{'})$ on the probability space $(W,{\mathcal W},Q_w)$.
We extend the domain of $Q^{(k)}$ to $W \times {\mathcal W}$
as follows:
\begin{equation} \label{2.14}
Q^{(1)} (w, B):= Q(w, B),\quad
Q^{(k)} (w, B):=
\int_{W} dQ_w(w^{'})\,Q^{(k-1)}(w', B), \quad k \geq 2.
\end{equation}
We will say $Q^{(k)}$ the \textit{$k$-step transition probability function
of the Markov chain} associated with $\{(W, {\mathcal W}), X, u, P\}$.
We see that $\left(W,{\mathcal W},Q^{(k)}_w\right)$ is also a probability space for each $w \in W$
and $k \geq 1$ where $Q^{(k)}_w:=Q^{(k)}(w,\cdot)$.
Furthermore, define
\begin{equation} \label{2.15}
Q_n(w,B):=\frac{1}{n}\sum_{k=1}^{n}Q^{(k)}(w,B),\quad
 (w,B)\in W \times {\mathcal W},\,n \geq 1.
\end{equation}
For example,
$Q_{1}=Q$.
Then $(W,{\mathcal W},Q_{n,w})$ is also a probability space
for each $w \in W$ where $Q_{n,w}(B):=Q_{n}(w,B)$ for $B \in {\mathcal W}$.
\subsection{Examples} \label{subsec2.3}
In this section, we shall give two examples of RSCC.
\begin{example} \label{ex2.4}
\rm
We show the RSCC associated with regular continued fraction expansions.
More precisely, it is the RSCC associated with the dynamical system
$(I,{\mathcal B}_I,\tau)$ where $\tau$ is the Gauss transformation, defined as $\tau = T_1$, where $T_N$ is as in (\ref{1.1}).
Define the RSCC $\{(W, {\mathcal W}), X,  u, P\}$ as
$W = [0,1], \quad {\mathcal W} = {\mathcal B}_{[0,1]},\quad
X = \mathbb{N}_+$,
\begin{equation} \label{2.16}
u: W \times X \rightarrow W; \quad u(w, x) = \frac{1}{w+x},
\end{equation}
\begin{equation} \label{2.17}
P : W \times X \rightarrow [0, 1]; \quad P(w, x) = \frac{w+1}{(w+x)(w+x+1)}.
\end{equation}
By definition, we see that $u(\cdot,x)$ is also a right inverse of $\tau$ for each $x \in X$, that is, $\tau(u(w,x))=w$ for any $w \in W$.
This shows that the dynamical system $(I,{\mathcal B}_I,\tau)$
is encoded onto the RSCC $\{(W, {\mathcal W}), X,  u, P\}$.
About more details, see Sec.1.2 of \cite{IG-2009}.
\end{example}
\begin{example} \label{ex2.5}
\rm
According to $\S$3.1 of \cite{Kawamura}, we show a simple (but non-trivial) example
of (finite) RSCC by using a finite automaton with input/output.
Define the data $\{(W, {\mathcal W}), X,  u, P\}$
as follows:
\begin{enumerate}
\item[(i)]
For the 2-point set $\{1,2\}$,
let ${\cal P}(\{1,2\})$ denote its power set.
Then $(\{1,2\}, {\cal P}(\{1,2\}))$ is a measurable space.
Define $(W, {\mathcal W}):=(\{1,2\}, {\cal P}(\{1,2\}))$
and $X:=\{1,2\}$.
\item[(ii)]
Any function $u: W \times X \rightarrow W$ is
measurable with respect to $(W,{\mathcal W})$ by definition.
For example,
let
\begin{equation} \label{2.18}
u(i,j):=j, \quad i,j=1,2.
\end{equation}
	\item[(iii)]
Any transition probability function $P:W \times X \to [0,1]$
is uniquely defined by two real numbers $P(i,1) \in [0,1]$ for $i=1,2$.
For example,
for $0\leq \alpha,\beta\leq 1$,
define
$P(1,1)=\alpha$ and $P(2,1)=\beta$.
Then they define a unique transition probability function $P$ as
\begin{equation} \label{2.19}
P(1,1)=\alpha,\quad
P(1,2)=1-\alpha,\quad
P(2,1)=\beta,\quad
P(2,2)=1-\beta.
\end{equation}
\end{enumerate}
For example,
the value $Q(1,1)$ of $Q$ in (\ref{2.9})
is computed as follows:
\begin{equation} \label{2.20}
Q(1,1)=P(1,\{x\in \{1,2\}:u(1,x)=1\})=P(1,1)=\alpha.
\end{equation}
Remark that $u$ in (\ref{2.18}) is the case such that
$u(i,\cdot)$ does not depend on $i \in W$ (\cite{IG-2009}, p15, (i)).
As the case such that $u(i,\cdot)$ depends on $i \in W$,
we can define other $u$, for example,
$u(i,j)=i$ for $i \in W$, $j \in X$.

We can illustrate this example as a finite automaton
with input/output (= Mealy machine \cite{Mealy}).
Assume that $\{(W, {\mathcal W}), X, u, P\}$
is as in (\ref{2.18}) and (\ref{2.19}).
For this system, we draw the transition diagram
(= an oriented graph with labeling edges) as follows:
\begin{enumerate}
\item[(1)]
The set of vertices is the state space $W = \{1,2\}$.
\item[(2)]
For two vertices $i,j \in W$ ($i$ and $j$ may be same),
if there exists an event (=input signal) $k \in X=\{1,2\}$ such that $u(i,k)=j$,
then draw the oriented edge from $i$ to $j$.
We write $(i,k,j)$ as this edge here.
\item[(3)]
Write ``$k/P(i,k)$" as the label of an edge $(i,k,j)$.
\end{enumerate}
According to these rules,
the transition diagram of the RSCC is illustrated as Figure \ref{figure:mealy}.

\noindent
%
\def\bcircle#1{\circle{300}\put(-30,-40){$#1$}}
\def\labee#1#2{\put(0,0){$#1/#2$}}
\def\numbers{\put(-800,0){\bcircle{1}}\put(800,0){\bcircle{2}}}
\def\edgethree{
\thicklines
\qbezier(-650,-30)(0,-200)(650,-30)
\qbezier(-650,30)(0,200)(650,30)
\put(-650,-20){\vector(-3,1){0}}
\put(670,20){\vector(3,-1){0}}
}
\def\edges{\thicklines\put(0,0){\edgethree}}
\def\semi{
\qbezier(150,0)(0,150)(-150,0)
\qbezier(150,0)(200,-75)(150,-150)
\qbezier(-150,0)(-200,-75)(-150,-150)
\put(-110,-170){\vector(1,-1){0}}}
\def\labels{
\put(-1320,-40){\labee{1}{\alpha}}
\put(1160,-40){\labee{2}{(1-\beta)}}
\put(-120,-270){\labee{1}{\beta}}
\put(-250,170){\labee{2}{(1-\alpha)}}
}
\def\graphs{
\thicklines
\put(0,0){\labels}
\put(0,0){\numbers}
\put(0,0){\edges}
%
\put(850,10){\rotatebox{-90}{\semi}}
\put(-1075,-20){\rotatebox{90}{\semi}}
}
%
%
\begin{fig} \label{figure:mealy}
\quad\\
\setlength{\unitlength}{.03mm}
\begin{picture}(1000,600)(0,500)
\put(2000,800){\graphs}
\end{picture}
\end{fig}

\noindent
For example,
from (\ref{2.18}),
$u(1,1)=1$, and
from (\ref{2.19}),
$P(1,1)=\alpha$.
Hence the label of the edge $(1,1,1)$ is given as
``$1/\alpha$".
About a Markov chain related to this example,
see Chap. 5 of \cite{Durrett}.
About other examples of Mealy machine,
see \cite{Kawamura,Kawamura2}.
\end{example}
\subsection{Operators associated with an RSCC} \label{subsec2.4}
In this subsection, we present the asymptotic and ergodic properties of
operators associated with an RSCC. To state these results, we prepare definitions.
We add the following assumption for all RSCC in this subsection:\\
\\
{\bf Assumption (B)}\,
For an RSCC $\{(W,{\mathcal W}),X,u,P\}$,
$W$ is a measurable subset of the measurable space
$(\mathbb {R},{\mathcal B}_{\mathbb{R}})$
and ${\mathcal W}={\mathcal B}_{W}$.\\

\noindent
Thanks to Assumption (B),
we can use the absolute value $|\cdot|$ and the Lebesgue measure
on $W$ induced by $\mathbb{R}$,
which will be necessary in this subsection.

Let $L^{\infty}(W)$ denote the Banach space of all
complex-valued bounded measurable functions defined on $(W,{\mathcal W})$.
We define operators on $L^{\infty}(W)$
associated with an RSCC $\{(W,{\mathcal W}),X,u,P\}$ as follows.
\begin{definition} \label{def2.7}
\begin{enumerate}
\item[(i)]
The \textit{transition operator} $U$ on $L^{\infty}(W)$ is defined by
\begin{equation} \label{2.21}
\{Uf\}(w) := \sum_{x \in X} P(w, x)\,f(u(w,x)), \quad
f \in L^{\infty}(W),\,w \in W.
\end{equation}
\item[(ii)]
(\cite{IG-2009}, (3.1.7))
For $n \geq 1$, define the operator $U_n$ on $L^{\infty}(W)$  as
\begin{equation} \label{2.22}
\{U_nf\}(w):= \int_{W} f(w')\, dQ_{n,w}(w^{'}),
\quad f \in L^{\infty}(W),\,w \in W
\end{equation}
where
$Q_{n,w}:=Q_{n}(w,\cdot)$ is as in (\ref{2.15}).
\end{enumerate}
\end{definition}
\begin{remark} \label{rem2.8}
{\rm
\begin{enumerate}
\item[(i)]
Let $U$ be as in (\ref{2.21}).
For each  $w \in W$,
$\{U(f)\}(w)$ is the expectation value of $(f\circ u)(w,\cdot)$
with respect to the probability space $(X,{\mathcal P}(X),P(w,\cdot))$.
For example, if $f$ is the characteristic function of $B \in {\mathcal W}$,
then we see $\{Uf\}(w)=Q(w,B)$.
For each $n\geq 1$,
$\{U^{n}f\}(w)=\sum_{x \in X^n}P_n(w,x)f(wx)$
where $U^n$ denotes the $n$-the iterate of $U$
and $P_n$ is as in (\ref{2.6}) for $r=n$.
\item[(ii)]
Let $U_n$ be as in (\ref{2.22}).
Remark that $U_n$ is well-defined on $L^{\infty}(W)$ because $Q_{n,w}$ is a probability measure.
For example, if $f$ is the characteristic function of $B\in {\mathcal W}$,
then $\{U_1 f\}(w)=Q(w,B)$.
Since
\begin{equation} \label{eqn2.23}
Q_{n,w}(w^{'})=Q_n(w,w^{'})=\frac{1}{n}\sum_{k=1}^{n}
Q^{(k)}(w,w^{'}),
\end{equation}
we see that
\begin{equation} \label{eqn2.24}
\{U_{n}f\}(w)=\frac{1}{n}
\sum_{k=1}^{n}\int_{W}f(w^{'}) \,dQ_{w}^{(k)}(w^{'}).
\end{equation}
\end{enumerate}
}
\end{remark}

Next, let  $L(W)$ denote the Banach space of all complex-valued Lipschitz continuous functions on $W$ with the following norm $\|\cdot\|_L$:
\begin{equation} \label{2.25}
\left\| f \right\|_L := \sup_{w \in W} |f(w)|
+ \sup_{w'\ne w''} \frac{|f(w') - f(w'')|}{|w' - w''|},
\quad f \in L(W).
\end{equation}
Remark that we use Assumption (B) here.
By definition, $L(W)\subset L^{\infty}(W)$.

According to \cite{IG-2009, Norman-1972}, we introduce several characterizations
of the operator $U$ in (\ref{2.21}) as follows.
\begin{definition}(\cite{IG-2009}, Definition 3.1.4, \cite{Norman-1972}, Definition 2.1) \label{def2.9}
Let $W, U,U_n, L(W)$ be as in (\ref{2.1}), (\ref{2.21}), (\ref{2.22}) and (\ref{2.25}), respectively.
We consider restrictions of $U,U_n$ on $L(W)$ as follows.
\begin{enumerate}
\item[(i)]
$U$ is ordered if there exists a bounded linear operator $S$ on $L(W)$ such that
\begin{equation} \label{2.26}
\lim_{n \rightarrow \infty} \|U_nf-S f\|_L = 0, \quad f \in L(W).
\end{equation}
\item[(ii)]
$U$ is aperiodic if there exists a bounded linear operator $S^{'}$ on $L(W)$ such that
\begin{equation} \label{2.27}
\lim_{n \rightarrow \infty} \|U^nf-S^{'} f \|_L = 0, \quad f \in L(W),
\end{equation}
where $U^n$ is the $n$-th iterate of $U$ for $n\geq 1$.
\item[(iii)]
$U$ is ergodic with respect to $L(W)$ if $U$ is ordered  and the rank of $S$ in (\ref{2.26}) is $1$.
\item[(iv)]
$U$ is regular with respect to $L(W)$ if $U$ is ergodic and aperiodic.
\item[(v)]
The Markov chain corresponding to $U$ is regular if $U$ is regular with respect to $L(W)$.
\end{enumerate}
\end{definition}
\begin{remark} \label{rem2.10}
{\rm
\begin{enumerate}
\item[(i)]
Definition \ref{def2.9}(i) and (ii) mean that sequences $(U_n)$
and $(U^n)$ converge to operators $S$ and $S^{'}$, respectively,
with respect to the strong operator topology on $L(W)$.
\item[(ii)]
Under the Assumption (B), $W$ is separable.
In addition, if $W$ is complete and $U$ is orderly,
then there exists probability measures $\{Q^{\infty}_w:w \in W\}$ on
$(W,{\mathcal W})$ such that
\begin{equation} \label{2.28}
\{Sf\}(w)=\int_W dQ^{\infty}_w(w^{'})\, f(w^{'}), \quad f\in L(W),\,w \in W
\end{equation}
by Theorem 3.1.24 of \cite{IG-2009}
where $S$ is as in Definition \ref{def2.9}(i).
\item[(iii)]
In addition to the assumption in (ii), if $U$ is ergodic with respect to $L(W)$,
then $Sf$ in (\ref{2.28}) is a constant function on $W$ for any $f \in L(W)$.
Therefore the operator $S$ is identified with a bounded linear functional
on $L(W)$, $S:L(W)\to \mathbb{C}$.
Then (\ref{2.28}) is rewritten as
\begin{equation} \label{2.29}
S:L(W)\rightarrow \mathbb{C};\quad
Sf=\int_W dQ^{\infty}(w^{'})\, f(w^{'}), \quad f \in L(W)
\end{equation}
for some probability measure $Q^{\infty}$ on $(W,{\mathcal W})$.
\end{enumerate}
}
\end{remark}

\begin{definition} \label{def2.11} (\cite{IG-2009}, Definition 3.1.15)
Under the Assumption (B), $\left\{(W,{\mathcal W}), X, u, P\right\}$ is
an RSCC with contraction if the following conditions are satisfied:
\begin{enumerate}
\item[(i)]
$r_1< \infty$,
\item[(ii)]
 $r_{\ell}<1$ for some $\ell \geq 1$, and
\item[(iii)]
$R< \infty$
\end{enumerate}
where $(r_k)$ and $R$ are defined as
\begin{eqnarray} \label{2.30}
r_k &:=& \sup_{w'\neq w''} \sum_{x \in X^k}P(w', x)
\frac{|w'x-w''x|}{|w'- w''|},\quad k\geq 1, \\
R &:=& \sup_{A \subset X}\, \sup_{w'\neq w''} \frac{|P(w', A)- P(w'', A)|}{|w'- w''|}. \label{2.31}
\end{eqnarray}
\end{definition}

\noindent
Remark that we use the Assumption (B) for $|\cdot|$ in
Definition \ref{def2.11}.

When $k=1$, (\ref{2.30}) is rewritten by using $u$ as follows:
\begin{equation} \label{2.32}
r_1 = \sup_{w'\neq w''} \sum_{x \in X}P(w', x)
\frac{|u(w',x)-u(w'',x)|}{|w'- w''|}.
\end{equation}
If $\sup_{w'\neq w''}|u(w',x)-u(w'',x)|/|w'- w''|<1$,
then $r_{1}<1$.
In this case, assumptions (i) and (ii) in
Definition \ref{def2.11} are satisfied.
\begin{theorem} \label{th2.12}
Under the Assumption (B), let $\left\{(W, {\mathcal W}), X, u, P\right\}$
be an RSCC with contraction. For $Q^{(n)}$ in (\ref{2.13}),
define $(\sigma_n)$ by
\begin{equation} \label{2.33}
\sigma_n(w) := \mathrm{supp}\, Q^{(n)}(w, \cdot),\quad w \in W
\end{equation}
where
$\mathrm{supp}\,\mu$ denotes the support of a measure $\mu$.
Assume that $W$ is compact.
Then the following holds.
\begin{enumerate}
\item[(i)]
The Markov chain associated with the RSCC is regular
if and only if there exists a point $w_0 \in W$ such that
\begin{equation} \label{2.34}
\lim_{n \rightarrow \infty} {\rm dist}(\sigma_n(w),w_0)=0\quad
\mbox{for all }w \in W
\end{equation}
where ${\rm dist}(A,w):=\inf_{w^{'}\in A}|w^{'}-w|$  for $A \subset W$.
\item[(ii)]
For $(\sigma_{n})$ in (\ref{2.33}),
the following holds:
\begin{equation} \label{2.35}
\sigma_{m+n}(w) = \overline{\bigcup_{w' \in \sigma_m(w)}\sigma_n(w')},
\end{equation}
for all $m$,
$n \in \mathbb{N}_+$, $w \in W$,
where the overline means the topological closure in $W$.
\end{enumerate}
\end{theorem}
%
%
\noindent \textbf{Proof.}
(i) See Theorem 3.3.31, p.116, of \cite{IG-2009}.

\noindent
(ii) By assumption, $W$ is a compact metric subspace of $\mathbb{R}$.
Hence the $Q$ in (\ref{2.10}) is ``continuous" in the sense of
Definition 3.3.1 in \cite{IG-2009}.
Therefore we can apply Lemma 3.3.32, p. 117 of \cite{IG-2009} and the statement holds.
\hfill $\Box$
\begin{definition} \label{def2.13}
Let $\{(W, {\mathcal W}), X, u, P\}$ be an RSCC and let
$P_r$ be as in (\ref{2.7}).
\begin{enumerate}
\item[(i)]
For $w \in W$, $n, r \in \mathbb{N}_+$, and $A \subset X^r$, define
\begin{equation} \label{2.36}
P^{n}_r(w, A) := P_{r+n-1}(w, \,X^{n-1} \times A),
\end{equation}
with the convention $X^{0} \times A := A$.
\item[(ii)]
(\cite{IG-2009}, Definition 2.1.4)
An RSCC $\left\{(W, {\mathcal W}), X, u, P\right\}$ is said to be uniformly ergodic if for any $r \in \mathbb{N}_+$, there exists a probability measure
$P_r^{\infty}$ on $(X^r,{\mathcal P}(X^r))$ such that
the sequence $\{P^{n}_r(w,A):n\geq 1\}$ in $(\ref{2.36})$ converges uniformly to  $P_r^{\infty}(A)$, that is,
the following sequence $(\varepsilon_n)_{n \in \mathbb{N}_+}$
goes to $0$ when $n \rightarrow \infty:$
\begin{equation} \label{2.37}
\varepsilon_n := \sup \{|P_r^{n}(w,A) - P_r^{\infty}(A)|
: w \in W,\, r \in \mathbb{N}_+,\, A \subset X^r \}.
\end{equation}
\end{enumerate}
\end{definition}
\begin{theorem} \label{th2.14}
Under the Assumption (B), let $\left\{(W,{\mathcal W}), X, u, P\right\}$ be an RSCC
with contraction such that $W$ is compact.
Assume that $\left\{(W, {\mathcal W}), X, u, P\right\}$
has a regular associated Markov chain.
Then the following holds:
\begin{enumerate}
\item[(i)]
$\left\{(W, {\mathcal W}), X, u, P\right\}$ is uniformly ergodic.
\item[(ii)]
Let $Q^{\infty}$ be as in Remark \ref{rem2.10}(iii).
Then the probability measure $P^{\infty}_r$ in  (\ref{2.37})
is given as follows:
\begin{equation} \label{2.38}
P^{\infty}_r(A) = \int_{W} P_r(w, A) \,dQ^{\infty}(w),
\quad A \in {\mathcal P}(X^{r})
\end{equation}
where $P_r$ is as in (\ref{2.7}).
\end{enumerate}
\end{theorem}

\noindent \textbf{Proof.}
(i) From Theorem 3.4.5, p.125 in \cite{IG-2009},
the statement holds.
\noindent
(ii) By (i), conditions in Remark \ref{2.10}(ii) and (iii) are satisfied by assumption. Hence we obtain $Q^{\infty}$ on $(W,\mathcal{W})$ in (\ref{2.29}).
From Theorem 3.4.5, p.125 in \cite{IG-2009},
the statement holds.
\hfill $\Box$\\
\section{Proof of Theorem \ref{G-K}} \label{section3}
In this section, we prove Theorem \ref{G-K}.
General results in Section \ref{section2} will be applied to $N$-continued fraction expansions.
\subsection{RSCC associated with $N$-continued fraction expansion} \label{subsec3.1}
Fix an integer $N \geq 1$.
We introduce a random system with complete connections (=RSCC)
$\{(I, {\mathcal{B}}_I), \Lambda, u, P\}$ as follows:
\begin{equation} \label{3.1}
\left\{
\begin{array}{ll}
u:I \times \Lambda \rightarrow I; \quad &u(x,i):=u_{N,i}(x)=\displaystyle \frac{N}{x+i},
\\
\\
P:I \times \Lambda \rightarrow I; \quad  &P(x,i):= V_{N,i}(x)=\displaystyle \frac{x+N}{(x+i)(x+i+1)},
\\
\\
\Lambda = \{N, N+1, \ldots \}.
\end{array}
\right.
\end{equation}

By the definition of $P$ and using the partial fraction decomposition,
$\sum_{i \in \Lambda} P(x,i)=1$.
By (\ref{2.10}),
$Q(x, B) = \sum_{i \in X(x,B)}V_{N,i}(x)$ for
$(x,B) \in I \times {\mathcal{B}}_I$
where $X(x,B) := \{i \in \Lambda: u_{N,i}(x) \in B \}$.
Let $U$ and $Q^{(n)}$ be as in (\ref{2.21}) and (\ref{2.14}), respectively.
By definition,
$\{(I, {\mathcal{B}}_I), \Lambda, u, P\}$ satisfies
the Assumption (B) in Section \ref{subsec2.4}.

For the dynamical system $(I,{\mathcal B}_I, T_N)$ in Section \ref{section1} and
a given probability measure $\mu$ on $(I,{\mathcal B}_I)$,
the ergodic behavior of the RSCC in (\ref{3.1})
allows us to find the limiting Gauss-Kuzmin distribution $F$
with respect to $(T_N, \mu)$:
\begin{eqnarray}
F(x) := \lim_{n \rightarrow \infty} \mu (T_N^n < x), \quad x \in I \label{3.2}
\end{eqnarray}
and the invariant measure induced by $F$.
\begin{lemma} \label{lem3.1}
$\{(I, {\mathcal{B}}_I), \Lambda, u, P\}$ is an RSCC with contraction.
\end{lemma}

\noindent \textbf{Proof.}
We have
\begin{eqnarray}\label{3.3}
  \frac{d}{dx}u(x,i) &=& \frac{d}{dx}u_{N,i}(x) = - \frac{N}{(x+i)^2}, \nonumber \\
   \frac{d}{dx}P(x,i) &=& \frac{i^2+i-2Ni - (x^2+2Nx+N)}{(x+i)^2(x+i+1)^2}
\end{eqnarray}
for any $x \in I$ and $i \in \Lambda$. Thus,
\begin{eqnarray}
  \sup_{x \in I}\left|\frac{d}{dx}u(x,i)\right| &\leq& \frac{N}{i^2}, \quad i \in \Lambda \\
  \sup_{x \in I}\left|\frac{d}{dx}P(x,i)\right| &<& \infty.
\end{eqnarray}
Hence the requirements of Definition \ref{def2.11} are fulfilled.

\hfill $\Box$

\begin{lemma} \label{lem3.2}
$\{(I, {\mathcal{B}}_I), \Lambda, u, P\}$ has a regular associated Markov chain.
\end{lemma}
\noindent \textbf{Proof.}
By Theorem \ref{th2.12}(i),
it is equivalent to that
\[(\Gamma):\quad
\left\{\begin{array}{l}
\mbox{there exists a point $x^* \in I$ such that}\\
\displaystyle\lim_{n \rightarrow \infty} {\rm dist}(\sigma_n(x),x^*)=0\quad
\mbox{for all }x \in I
\end{array}
\right\}
\]
where we remark that $W=I=[0,1]$ in this case.
Hence we show the condition $(\Gamma)$ as follows.

Fix $x \in I$.
Let us define
the sequence $(x_n)_{n \geq 0}$ in $I$,
recursively by $x_0:=x$, $x_{n+1} := \displaystyle \frac{N}{x_n + N}$ for
$n \geq 1$. Clearly $x_{n+1} \in \sigma_1(x_n)$
and therefore Theorem \ref{th2.12}(ii) and an induction argument lead us to the conclusion that $x_n \in \sigma_n(x)$
for $n \in \mathbb{N}_+$.
But, $\displaystyle \lim_{n \rightarrow \infty} x_n = x^*= \frac{-N+\sqrt{N^2+4N}}{2}$
for any $x \in I$. Hence
 ${\rm dist}(\sigma_n(x),\,x^*) \leq
|x_n - x^*| \rightarrow 0$ as  $n \rightarrow \infty$.
Hence we find $x^*:=\displaystyle \frac{-N+\sqrt{N^2+4N}}{2}$ in the condition $(\Gamma)$.
\hfill $\Box$

\begin{proposition} \label{prop3.3}
$\{(I, {\mathcal{B}}_I), \Lambda, u, P\}$ is uniformly ergodic.
\end{proposition}
\noindent \textbf{Proof.}
In order to apply Theorem \ref{th2.14}
to the RSCC $\{(I, {\mathcal{B}}_I), \Lambda, u, P\}$
in (\ref{3.1}), we check assumptions in Theorem \ref{th2.14}.
By definition, $I$ is compact. By Lemma \ref{lem3.1},
$\{(I, {\mathcal{B}}_I), \Lambda, u, P\}$ is an RSCC with contraction.
By Lemma \ref{lem3.2}, $\{(I, {\mathcal{B}}_I), \Lambda, u, P\}$ is regular.
Hence all assumptions in Theorem \ref{th2.14} are satisfied.
Hence the statement holds.

\hfill $\Box$

Let $L(I)$ denote the Banach space of all complex-valued Lipschitz continuous functions on $I$.
The regularity of $U$ in (\ref{1.4}) with respect to $L(I)$ follows from Theorem \ref{th2.12}.
Moreover, the sequence $\{Q^{(n)} (\cdot, \cdot):n \geq 1\}$
in (\ref{2.14}) converges uniformly to a probability measure $Q^{(\infty)}$ on $(I, {\mathcal B}_I)$ and that there exist two positive constants $q < 1$ and $k$ such that
\begin{equation}
\|U_n f-U_{\infty} f \|_L \leq k q^n \|f \|_L, \quad
n \in \mathbb{N}_+,\, f \in L(I) \label{3.6}
\end{equation}
where
\begin{eqnarray} \label{3.7}
U_n: L(I) \rightarrow L(I); \quad
\{U_n f\}(x) &=& \int_{I}f(y)\,dQ^{(n)}_x(y), \\
U_{\infty}:L(I) \rightarrow \mathbb{C}; \quad
U_{\infty} f &=& \int_{I}f(y)\,dQ^{(\infty)}(y) \label{3.8}
\end{eqnarray}
with $Q^{(n)}_x(B):=Q^{(n)}(x,B)$ in (\ref{2.14}).
\begin{proposition} \label{prop3.4}
The probability $Q^{(\infty)}$ is the invariant probability measure of the transformation $T_N$.
\end{proposition}
\noindent \textbf{Proof.}
For $G_N$ in Section \ref{section1} and $Q$ in (\ref{2.10}), and on account of the uniqueness of $Q^{(\infty)}$ we have to show that
\begin{equation} \label{3.9}
\int_{0}^{1} Q(x, B) \,dG_N(x) = G_N(B), \quad B \in {\mathcal B}_I.
\end{equation}
Since the intervals $\left[0, u\right) \subset \left[0, 1\right)$
generate ${\mathcal B}_I$, it is sufficient to show the equation (\ref{3.9})
just for $B = \left[0, u\right)$, $0 < u \leq 1$. Let
$E(x,N) = \left\lfloor\frac{N}{u}-x\right\rfloor +1$.
Since $Q(x, B) = \sum_{i \in X(x,B)}V_{N,i}(x)$ for $(x,B)
\in I \times {\mathcal{B}}_I$
where $X(x,B) := \{i \in \mathbb{N}: u_{N,i}(x) \in B \}$, we have
\begin{eqnarray} \label{3.11}
Q_m(x, [0, u)) &=& \sum_{\left\{i \in \mathbb{N}:
0 \leq u_i(x) < u \right\}}V_{N,i}(x) = \sum_{i \geq E(x,N)}V_{N,i}(x) \nonumber \\
&=& \frac{N-E(x,N)}{x+E(x,N)}.
\end{eqnarray}
Thus,
\begin{equation} \label{3.12}
\int_{0}^{1} Q (x, [0, u)) dG_N(x) = \frac{1}{\log \{(N+1)/{N}\}} \log \frac{x+N}{N}= G_N([0, u)).
\end{equation}
Hence the statement holds.
\hfill $\Box$
\subsection{Proof of Theorem \ref{G-K}} \label{subsec3.2}
By using Proposition \ref{prop3.3}, we prove Theorem \ref{G-K}.\\

\noindent \textbf{Proof of Theorem \ref{G-K}\,}
Fix an integer $N \geq 1$ and let $G_N$ be as in Section \ref{section1}.
By (\ref{1.5}), we have
\begin{equation} \label{3.13}
\mu \left(T_N^{-n}(A)\right) = \int_{A} \{U^nf_0\}(x)dG_N(x) \quad
\mbox{for any }n \in \mathbb{N}_+,\, A \in {\mathcal{B}}_I
\end{equation}
where $f_0(x)=(x+N)(d\mu / d\lambda) (x)$ for $x \in I$.
If $d\mu / d\lambda \in L(I)$, then $f_0 \in L(I)$
and by (\ref{3.8}) we have
\begin{eqnarray}
U^{(\infty)} f_0 &=& \int_{I} f_0(x)\,Q^{(\infty)}(dx) = \int_{I} f_0(x)\,G_N(dx) \nonumber \\
&=& \int_{I} (d\mu / d\lambda) (x)\,dx = \frac{1}{\log \{(N+1)/{N}\}} \log \frac{x+N}{N}. \label{3.14}
\end{eqnarray}

Taking into account (\ref{3.6}), there exist two constants $q<1$ and $k$ such that
\begin{equation} \label{3.15}
\left\|U^n f_0 - U^{(\infty)} f_0\right\|_L \leq k q^n\left\|f_0\right\|_L,
\quad n \in \mathbb{N}_+.
\end{equation}
Furthermore, consider the Banach space $C(I)$ of all real-valued continuous functions on $I$ with the norm
$\|f \| := \sup_{x \in I}|f(x)|$.
Since $L(I)$ is a dense subspace of $C(I)$ we have
\begin{equation} \label{3.16}
\lim_{n \rightarrow \infty} \left\|\left(U^n - U^{(\infty)}\right)f\right\| = 0 \quad \mbox{for all } f \in C(I).
\end{equation}
Therefore, (\ref{3.16}) is valid for a measurable function $f_0$ which is $Q^{(\infty)}$-almost surely continuous, that is, for a Riemann-integrable function $f$.
Thus, we have
\begin{eqnarray} \label{3.17}
\lim_{n \rightarrow \infty} \mu \left(G^n_N < x\right)
&=& \lim_{n \rightarrow \infty} \int_{0}^{x} \{U^nf_0\}(u) G_N(du) \\
\label{3.18}
&=& \frac{1}{\log \{(N+1)/{N}\}} \log \frac{x+N}{N} \int_{0}^{x} \,G_N(du)\\
\label{3.19}
&=& \frac{1}{\log \{(N+1)/{N}\}} \log \frac{x+N}{N}.
\end{eqnarray}
Hence (\ref{1.6}) is proved.
\hfill $\Box$
\\

\end{document}